\documentclass[preprint,12pt]{elsarticle}
\usepackage{etoolbox}
\makeatletter
\def\ps@pprintTitle{%
	\let\@oddhead\@empty
	\let\@evenhead\@empty
	\def\@oddfoot{\reset@font\hfil\thepage\hfil}
	\let\@evenfoot\@oddfoot
}
\makeatother
\usepackage{graphics}
\usepackage{amsmath, amsthm, amssymb}
\usepackage{natbib}
\usepackage[colorlinks=true]{hyperref}
\usepackage{lscape}
\usepackage{tabularx}
\usepackage[utf8]{inputenc} 
\usepackage{amsthm}

\usepackage{enumerate}
\usepackage{mathrsfs}
\usepackage{setspace}
\usepackage{multirow}
\usepackage{graphicx}
\usepackage{amsfonts}
\usepackage{verbatim}
\usepackage{amssymb}
\usepackage{amsthm}
\usepackage{caption}
\usepackage{subcaption}
\usepackage{multirow,bigstrut,bigdelim}

\theoremstyle{plain}

\newtheorem{ex}{Example}[section]

\newtheorem{prop}{Proposition}[section]
\theoremstyle{remark}
\newtheorem{rem}{Remark}[section]

\numberwithin{equation}{section}
\usepackage[mathscr]{euscript}
\usepackage{geometry} 
\usepackage{graphicx,lscape} 
\usepackage{booktabs}
\usepackage{array}
\usepackage{paralist} 
\usepackage{verbatim}
\usepackage{lscape}
\usepackage{makecell}
\biboptions{comma,round,authoryear}

\begin{document}
	\begin{frontmatter}
		\title{\textbf{The Relative Information Generating Function\\
        -A Quantile Approach}}		\author{P. G. Sankaran*\footnote{*\it Corresponding author}, S. M. Sunoj,  and Pavithra Hariharan }
       \ead{sankaran.p.g@gmail.com, smsunoj@cusat.ac.in, pavithrahariharan97@gmail.com}
		\address{Department of Statistics\\Cochin University of Science and Technology\\ Kerala, India 682 022.}
%         A Quantile Approach}}		
%         \author{
% \name{P. G. Sankaran\thanks{CONTACT P. G. Sankaran. Email: sankaran.pg@gmail.com} , S. M. Sunoj and Pavithra Hariharan }
% \affil{Department of Statistics, Cochin University of Science and Technology, Cochin 682 022, Kerala, India }
% }

		\begin{abstract}
         Information generating functions have been used for generating various entropy and divergence measures. In the present work, we introduce quantile based relative information generating function and study its properties. The proposed generating function provides well-known Kullback-Leibler divergence measure. The quantile based relative information generating function for residual and past lifetimes are presented. A non parametric estimator for the function is derived. A simulation study is conducted to assess performance of the estimators. Finally, the proposed method is applied to a real life data.
		\end{abstract}
	\end{frontmatter}
	\section{Introduction}
Divergence measures play a pivotal role to measure the closeness between two probability distributions. One of the best known measures is the Kullback-Leibler (K-L) divergence measure which quantify the discrepancy between two probability distributions (\cite{kullback1951information}). Many generalizations of these measures have been studied in literature by several authors (see \cite{ali1966general}, \cite{csiszar1967two}, \cite{lin1991divergence}, \cite{di2015some}, \cite{kharazmi2018time}). The applications of these measures include in the analysis of contingency tables (\cite{gokhale1978minimum}), in approximation of probability distributions (\cite{chow1968approximating} and \cite{kazakos1980decision}), in signal processing (\cite{kadota1967best} and \cite{kailath1967divergence}) and in pattern recognition (\cite{ben1978irrelevant},\cite{chen1973statistical}, and \cite{chen1976statistical}).

Like moment generating function in distribution theory, information generating functions have been defined for probability densities to determine entropy measures and divergence measures. \cite{golomb1966run} proposed an information generating function of probability density functions, whose derivative evaluated at one, gives well known Shannon's entropy measure. Relative entropy, based on Golomb information generating function has been defined by \cite{guiasu1985relative}, and its first derivative at one leads to K-L divergence. Let $X_1$ and $X_2$ be two random variables with probability densities $f_1(\cdot)$ and $f_2(\cdot)$ respectively. The relative information generating function, for any $\alpha>0$ is defined as 
\begin{equation}\label{1.1}
I_\alpha(f_1,f_2)=\int_{}^{}f_1(x)^{\alpha}f_2(x)^{1-\alpha}dx,
\end{equation}
provided the integral exists. 
Note that for $\alpha=1$, $I_1(f_1,f_2)=1$ and the derivative of \eqref{1.1} with respect to $\alpha$ at one is
\begin{equation}\label{1.2}
    \frac{d }{d\alpha}I_\alpha(f_1,f_2)|_{\alpha=1}=\int_{}^{}f_1(x)\log\frac{f_1(x)}{f_2(x)}dx,
\end{equation}
which is well known K-L divergence measure. For more properties and applications on divergence measure, one could refer to \cite{lin1991divergence}, \cite{kharazmi2021cumulative}, and \cite{kharazmi2024jensen}.

The measure \eqref{1.1} has been developed and studied by using the properties of 
probability density functions, but in many situations \eqref{1.1} is not suitable as distribution functions (density functions) are not analytically tractable. An alternative method is to use quantile function of $X_i$ defined by 
\begin{equation}
Q_i(u)=F_i^{-1}(u)=inf\left\{x|F_i(x)\geq u\right\};\quad 0\leq u \leq 1,
\end{equation}
where $F_i(x)$ is the distribution function of $X_i$, $i=1,2$.

It is shown that quantile functions are less influenced by extreme observations so that the quantile based approach provides robust inference with a limited amount of information. Further, quantile functions have certain unique properties which are useful for modeling and analysis of data. For instance, sum of two distribution functions need not be a distribution function, but sum of two quantile functions is always a quantile function. \cite{gilchrist2000statistical} and \cite{nair2013quantile} provided extensive review on quantile based analysis of statistical data.

Motivated by this \cite{sankaran2016kullback} introduced quantile version of K-L divergence and studied its properties. \cite{sunoj2018quantile} discussed cumulative divergence measure in the quantile setup. \cite{sunoj2024quantile} presented quantile based cumulative K-L divergence in past time and discussed various applications. \cite{mansourvar2022phi} proposed a $\phi$ divergence family of measures in the quantile setup and studied its properties. In the present work, we introduce quantile version of relative information generating function and study its properties. The proposed function represents a family of divergence measures, including the well-known K-L divergence measure.

The rest of the paper is organized as follows. In Section \ref{sec2}, we define relative information generating function in the quantile setup. Various properties of the function are also discussed. Section \ref{sec3} develops the quantile form of the relative information generating function for residual lifetime distributions.
The quantile based relative information generating function for the past lifetime distribution is discussed in Section \ref{sec4}. A non parametric estimator for the relative information generating function is derived in Section \ref{sec5}. Finally, Section \ref{sec6} assesses performance of the estimators through a simulation study and Section \ref{sec7} illustartes its application to prostate cancer data.

\section{Relative Information Generating Function}\label{sec2}
Let $X_1$ and $X_2$ be two non-negative random variables with absolutely continuous distribution functions $F_1(x)$ and $F_2(x)$ and corresponding probability
densities $f_1(x)$ and $f_2(x)$. Denote $Q_i(x)$ as quantile functions of $X_i, i=1,2$. Note that $F_i\left(Q_i(p)\right)=p$ and $q_i(p) f_i\left(Q_1(p)\right)=1$, where $q_i(p)$ is the quantile density function of $X_i, i=1,2$. Setting $x=Q_1(p)$ in \eqref{1.1}, the quantile version of the relative information generating function is defined by
%\begin{equation}\label{2.1}
\begin{align}
I_\alpha^*\left(q_1, q_2\right) & =\int_0^1\left(\frac{f_1\left(Q_1(p)\right)}{f_2\left(Q_1(p)\right)}\right)^\alpha f_2\left(Q_1(p)\right) dQ_1(p) \\\label{2.1}
= & \int_0^1\left[\frac{1}{q_1(p) f_2\left(Q_1(p)\right)}\right]^{\alpha}f_2(Q_1(p)) q_1(p) dp\\\label{2.2}\nonumber
\end{align}
%\end{equation}
Since $F_1^{-1}(p)=Q_1(p)$, we obtain $F_2\left(F_1^{-1}(p)\right)=F_2\left(Q_1(p)\right)=Q_2^{-1}\left(Q_1(p)\right)$. The differentiation of the above with respect to $p$ implies that
\begin{equation*}
\frac{d}{d p}\left(F_2\left(F_1^{-1}(p)\right)\right)=\frac{d}{d p}\left(F_2\left(Q_1(p)\right)\right. =q_1(p) f_2(Q_1(p)) =\frac{d}{d p}\left(Q_2^{-1}\left(Q_1(p)\right)\right).
\end{equation*}
Denote $Q_3(p)=Q_2^{-1}\left(Q_1(p)\right)$. Then we have $q_3(p)=\frac{d Q_3(p)}{d p}=q_1(p) f_2\left(Q_1(p)\right)$. Now \eqref{2.1} becomes
\begin{equation}\label{2.2}
I_\alpha^*=I_\alpha^*\left(q_1, q_2\right)=\int_0^1 (q_3(p))^{1-\alpha}dp, \quad \alpha>0.    
\end{equation}
It may be noted that $q_3(p)$ is the quantile density of $F_2\left(Q_1(p)\right)=Q_2^{-1}\left(Q_1(p)\right)$. The measure $I_\alpha^{*}$ is always non-negative.

\begin{rem}
It is easy to show that Kullback-Leibler (K-L) divergence measure is a special case of \eqref{2.2}. Differentiating \eqref{2.2} with respect to $\alpha$, we get 
\begin{equation}
  \frac{\partial I_\alpha^*}{\partial \alpha}=-\int_0^1 \log q_{3}(p) (q_3(p))^{1-\alpha} dp.
\end{equation}
Then
\begin{equation}\label{KL}
\left.\frac{\partial I_\alpha^*}{\partial \alpha}\right|_{\alpha=1}=-\int_0^1 \log q_3(p) dp,    
\end{equation}
which is well known K-L divergence measure in the quantile setup, introduced by \cite{sankaran2016kullback}.
\end{rem}

\begin{rem}
Differentiating \eqref{2.2} with respect to $\alpha,k$ times and substituting $\alpha=1$, we get
\begin{equation}
 \left.\frac{\partial^k I_\alpha^*}{\partial \alpha^k}\right|_{\alpha=1}=\int_0^1\left(-\log q_3(p)\right)^k dp,   
\end{equation}
which is the quantile version of generalized K-L divergence.
\end{rem}

\begin{rem}
Let $H_i(p)=\frac{1}{(1-u) q_i(p)}$ be the hazard quantile functions of $X_i,  i=1,2$ which is the quantile based hazard rates introduced in \cite{nair2009quantile}. Then we can write
\begin{equation}
 I_\alpha^*=\int_0^1\left[(1-p) H_3(p)\right]^{\alpha-1} dp,   
\end{equation}
where $H_3(p)$ is the hazard quantile function corresponding to $Q_3(p)$. On similar lines, one can write
\begin{equation}
I_\alpha^*=\int_0^1\left[p \tilde{H}_3(p)\right]^{\alpha-1} d p,
\end{equation}
where $\tilde{H}_3(p)$, is the reversed hazard quartile function corresponding to $Q_3(p)$. This is the quantile version of the reversed hazard rate, introduced in \cite{nair2013quantile}.
\end{rem}

\begin{rem}
Let $\alpha=\eta+\beta.$ We can write 
\begin{equation*}
    I_\alpha^*=e^{(\eta-\beta)D_\eta^\beta(X_1,X_2)},
\end{equation*}
where $D_\eta^\beta(X_1,X_2)$ is the quantile divergence measure of order $(\eta,\beta)$ introduced by \cite{kumar2018quantile}.
\end{rem}

\begin{rem}
  The Hellinger distance measure between two probability density functions $f_1(x)$ and $f_2(x)$ (\cite{hellinger1909neue}) is given by 
\[
H(f_1, f_2) = \frac{1}{2} \int _{0}^{\infty}\left(\sqrt{f_1(x)} - \sqrt{f_2(x)}\right)^2 dx
\]
\[
= 1 - \int _{0}^{\infty}\sqrt{f_1(x) f_2(x)} \, dx.
\]
In the quantile setup, by putting $x = Q(p)$, we get
\[
H^*(q_1, q_2) = 1 - \int_{0}^{1}  (q_3(p))^{1/2} \, dp.
\]
It may be noted that 
\[
H^*(q_1, q_2) = 1 - I^*_\frac{1}{2}.
\]
\end{rem}

\begin{rem}
The Bhattacharyya distance between $f_1(x)$ and $f_2(x)$ (\cite{kailath1967divergence}) is defined by
\[
B(f_1, f_2) = -\log \int _{0}^{\infty}\sqrt{f_1(x) f_2(x)} \, dx.
\]
The quantile version of $B(f_1, f_2)$ is obtained as
\[
B^*(q_1, q_2) = -\log \int _{0}^{1}  (q_3(p))^{1/2} \, dp,
\]
which can be written as
\[
B^*(q_1, q_2) = -\log I^*_\frac{1}{2}.
\]
\end{rem}

\begin{rem}
    The function \eqref{2.2} can be related to other divergence measures such as R$\acute{e}$nyi's divergence measure (\cite{renyi1961measures}), Tsallis divergence and alpha divergence measures (\cite{havrda1967quantification}), Sharma-Mittal divergence measure (\cite{sharma1975new}), and Cressie-Read divergence measure (\cite{cressie1984multinomial}). For instance, the well known R$\acute{e}$nyi's information measure of order $\alpha (\alpha\neq 1)$ between $f_1(x)$ and $f_2(x)$ is defined by
    \begin{equation*}
        D_\alpha(f_1,f_2)=\frac{1}{\alpha-1}\log\int_{0}^{\infty}f_1(x)^\alpha
        f_2(x)^{1-\alpha}dx
    \end{equation*}
(\cite{van2014renyi}). The quantile version of $D_\alpha(f_1,f_2)$ is
\begin{equation*}
    D_\alpha^*(q_1,q_2)=\frac{1}{\alpha-1}\log\int_{0}^{1}q_3(p)^{1-\alpha}dp.
\end{equation*}
Now, $I_\alpha^*$ is related to $D_\alpha^*(q_1,q_2)$ by
\begin{equation*}
        D_\alpha^*(q_1,q_2)=\frac{1}{\alpha-1}\log \,I_\alpha^*.
\end{equation*}
\end{rem}

\begin{prop}
  Suppose that $X_i$ has quantile density $q_i(p);i=1,2$. Let $Q_3(p)=Q_2^{-1}(Q_1(p))$ and $q_3(p)=\frac{dQ_3(p)}{dp}$. Then
  \begin{equation*}
I_\alpha^*=\sum_{k=0}^{\infty}\frac{(1-\alpha)^k}{k!}S_k(q_3),
  \end{equation*}
  where $S_k(q_3)=\int_{0}^{1}(\log q_3(p))^k dp$.
  \begin{proof}
      We can write
      \begin{equation*}
          I_\alpha^*=\int_{0}^{1}e^{(1-\alpha)\log q_3(p)}dp.
      \end{equation*}
    Using Maclaurin's theorem, we obtain
    \begin{equation*}
        \begin{aligned}
          I_\alpha^*&= \int_{0}^{1}\sum_{k=0}^{\infty}\frac{(1-\alpha)^k}{k!}(\log q_3(p))^k dp\\
          &=\sum_{k=0}^{\infty}\frac{(1-\alpha)^k}{k!}\int_{0}^{1}(\log q_3(p))^k dp\\
          &=\sum_{k=0}^{\infty}\frac{(1-\alpha)^k}{k!}S_k(q_3).
        \end{aligned}
    \end{equation*}
  \end{proof}
\end{prop}
We now present lower and upper bounds for $ I_\alpha^*$ in terms of K-L divergence and hazard quantile functions.
\begin{prop}
    Suppose that  $I_\alpha^*$ is the relative information generating function between $q_1(\cdot)$ and $q_2(\cdot)$. Then we have 
    \begin{equation*}
        L(\alpha)\leq I_\alpha^*\leq  U(\alpha),
    \end{equation*}
    where $L(\alpha)=\max(0,(1-\alpha)S_1(q_3))$ and $U(\alpha)=\int_{0}^{1}(H_3(p))^{\alpha-1}dp$, where $S_1(q_3)$ is quantile K-L divergence. 
    \begin{proof}
        We use the inequality $x^{1-\alpha}\geq (1-\alpha)\log x+1$. Substituting this inequality in \eqref{2.2}, we get the lower bound. The upper bound can be obtained using the relations $1-p\leq 1$ and $q_3(p)=\frac{1}{(1-p)H_3(p)}$.
    \end{proof}
\end{prop}
\begin{ex}
    Suppose that $X_i$ has exponential distribution with mean parameter $\lambda_i, i=1,2$. Then
$Q_i(p)=-\lambda_i\log (1-p), \lambda_i>0, i=1,2$. Then
$Q_3(p)=Q_2^{-1}\left(Q_1(p)\right)=1-(1-p)^{\lambda_1 /\lambda_ 2}$. Then  $q_3(p)=\frac{\lambda_1}{\lambda_2}(1-p)^{\frac{\lambda_1}{\lambda_2}-1}$.
We get 
\begin{equation*}
    \begin{aligned}
I_\alpha^*&=\int_0^1\left(\frac{\lambda_1}{\lambda_2}\right)^{1-\alpha}(1-p)^{(\frac{\lambda_1}{\lambda_2}-1)(1-\alpha)}dp\\
%&=\left(\frac{\lambda_1}{\lambda_2}\right)^{1-\alpha} \frac{(1-p)^{(\frac{\lambda_1}{\lambda_2}-1)(1-\alpha)}+1}{-\left[\left(\frac{\lambda_1}{\lambda_2}-1\right)(1-\alpha)+1\right]} \Bigg|_{0}^{1}\\
%&=\left(\frac{\lambda_1}{\lambda_2}\right)^{1-\alpha} \frac{\lambda_1}{\left(\lambda_1-\lambda_2\right)(1-\alpha)}+\lambda_2\\
&=\left(\frac{\lambda_1}{\lambda_2}\right)^{1-\alpha}\frac{\lambda_2}{\lambda_1+\alpha(\lambda_2-\lambda_1)}.
\end{aligned}
\end{equation*}
\end{ex}

In lifetime studies, Cox proportional hazards model plays an important role to study the effect of covariates on lifetime. The model is defined by $\overline{F}_2(x)=(\overline{F}_1(x))^\theta;\theta>0$, where $\overline{F}_1(\cdot)$ and $\overline{F}_2(\cdot)$ are survival functions of $X_1$ and $X_2$. The quantile version of the model is given by $Q_2(p)=Q_1(1-(1-p)^{1/\theta}); \theta>0$ (\cite{nair2018proportional}).

\begin{ex}
Suppose that $X_1$ and $X_2$ satisfy the property of Cox proportional hazards model. Then 
$Q_3(p)=1-(1-p)^{\theta}$ and $q_3(p)=\theta(1-p)^{\theta-1}$. Now
\begin{equation*}
\begin{aligned}
    I_\alpha^*&=\int_{0}^{1}\theta^{1-\alpha}(1-p)^{(\theta-1)(1-\alpha)}dp\\
    %&=\theta^{1-\alpha}\frac{(1-p)^{(\theta-1)(1-\alpha)}+1}{-[(\theta-1)(1-\alpha)+1]}\Biggr|_{0}^{1}\\
    %&=\theta^{1-\alpha}\frac{1}{\theta-1-\alpha\theta+\alpha+1}\\
    &=\frac{\theta^{1-\alpha}}{\alpha+\theta(1-\alpha)}.
\end{aligned}
\end{equation*}   
\end{ex}

\begin{ex}
 Let $X_1$ follows power-Pareto distribution defined by $Q_1(p)=c p^{\lambda_1}(1-p)^{-\lambda_2},\quad c, \lambda_1, \lambda_2>0$ and $X_2$ follows power distribution with $Q_2(p)=\beta_1 p^{1 / \beta_2}, \quad \beta_1, \beta_2>0$. Then we obtain $q_3(p)=\frac{\beta_2 c^{
\beta_2}}{\beta_1^{\beta_2}} p^{\lambda_1 \beta_2-1}(1-p)^{\lambda_2 \beta_2-1}\left[\lambda_1+p(\lambda_2-\lambda_1)\right]$. Then,
\begin{equation*}
   I_\alpha^{*}=\int_{0}^{1}  \beta_2^{1-\alpha}\left(\frac{c}{\beta_1}\right)^{\beta_2(1-\alpha)}
   p^{(\lambda_1 \beta_2-1)(1-\alpha)}(1-p)^{(\lambda_2 \beta_2-1)(1-\alpha)}\left[\lambda_1+p(\lambda_2-\lambda_1)\right]^{1-\alpha}dp,
\end{equation*}
which can be evaluated numerically for specific values of $ \beta_1, \beta_2,c, \lambda_1,$ and $ \lambda_2$. Plot of $I_\alpha^{*}$ versus $\alpha$ for various combinations of $(\beta_1,\beta_2,c,\lambda_1,\lambda_2)$ is shown in Figure \ref{Ialpha}. It is clear that $I_\alpha^{*}$ increases as $\alpha$ increases.

\begin{figure}[h!]
\centering
\includegraphics[width=8cm,height=5cm]{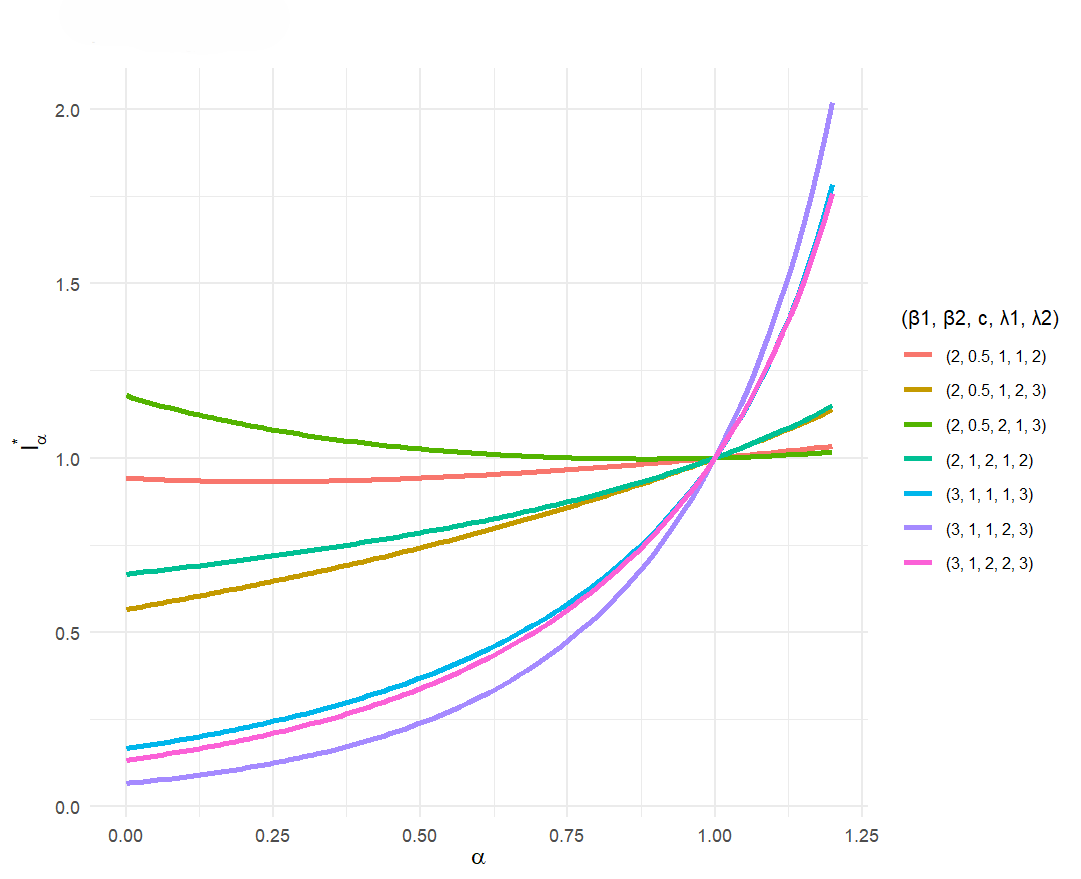}
\caption{Plot of $I_{\alpha}^*$ versus $\alpha$ for various $(\beta_1,\beta_2,c,\lambda_1,\lambda_2)$}
\label{Ialpha}
\end{figure}

\end{ex}

One of the objectives of lifetime studies is to identify the best quantile function for the given set of data. This can be done by choosing a basic quantile function like uniform, exponential etc. and applying suitable transformation to arrive at a flexible model. We now present the relative information generating function \eqref{2.2} in the context of transformation of quantile functions.

Let $T_1(\cdot)$ and $T_2(\cdot)$ be two continuous, non-decreasing,  and invertible transformations. Then the relative information generating function for the probability distributions of $T_1(X_1)$ and $T_2(X_2)$ is given by
\begin{equation*}
 I_\alpha^*\left(T_1\left(X_1\right), T_2\left(X_2\right)\right)=\int_0^1\left[\frac { d } { d p } \left[Q_2^{-1}\left[T_2^{-1}\left(T_1\left(Q_1\left(p\right)\right)\right)\right]\right]\right]^{1-\alpha} dp.  
\end{equation*}
\begin{ex}
 Let $X_1$ and $X_2$ follow two independent Pareto I distributions with quantile functions $Q _i(p)=(1-p)^{-\gamma_i}, \quad \gamma_i>0, i=1,2.$ Take $T_1\left(X_1\right)=\log X_1$ and $T_2\left(X_2\right)=\log X_2$. Then the relative information generating function for two exponential distributions with means $\gamma_1$ and $\gamma_2$ can be obtained as
 \begin{equation*}
 %\begin{aligned}
I_\alpha^*(T_1(\cdot),T_2(\cdot))=\int_{0}^{1}\left(\frac{\gamma_1}{\gamma_2}(1-p)^{\frac{\gamma_1}{\gamma_2}-1}\right)^{1-\alpha}dp=\left(\frac{\gamma_1}{\gamma_2}\right)^{1-\alpha}\frac{\gamma_2}{(\gamma_2-\gamma_1)\alpha+\gamma_1}.
%\end{aligned}    
\end{equation*}
\end{ex}

\section{Relative Information Generating Function for Residual Lifetimes}\label{sec3}
In reliability analysis, the experimenter has information about the current age $t$ of the systems. In such context, $I_\alpha^*(\cdot)$ is not an appropriate function for studying the divergence between the probability distributions. Suppose that $X_1$ and $X_2$ have truncated at time $t$. Then the residual lifetimes are  $X_1^t=X_1-t \mid X_1 >t$ and $X_2^t=X_2-t \mid X_2>t$. The residual lifetime survival distributions of $X_1^t$ and $X_2^t$ are
$$
\overline{F}_1^t(x)=\frac{\overline{F}_1(x+t)}{\overline{F}_1(t)}; \quad x>t
$$
and
$$
\overline{F}_2^t(x)=\frac{\overline{F}_2(x+t)}{\overline{F}_2(t)}; \quad x>t
$$
where $\overline{F}_i(\cdot)$ is the survival function of $X_i, i=1,2$. Then the probability density function of $X_i^t$ is
$$
\begin{aligned}
f_i^t(x)=\frac{f_i(x+t)}{\overline{F}_i(t)}, \quad x>t>0, \quad i=1,2.
\end{aligned}
$$
The relative information generating function between $X_1^t$ and $X_2^t$ is given by
\begin{equation}\label{eq3.1}
R_\alpha(t)  =\int_0^{\infty}\left[\frac{f_1(x+t) /\overline{F}_1(t)}{f_2(x+t) / \overline{F}_2(t)}\right]^\alpha \frac{f_2(x+t)}{\overline{F}_2(t)} d x 
 =\int_t^{\infty}\left(\frac{f_1(y)}{f_2(y)} \frac{\overline{F}_2(t)}{\overline{F}_1(t)}\right) ^\alpha\frac{f_2(y)}{\overline{F}_2(t)} d y
\end{equation}
Let $F_i\left(Q_i(p)\right)=p,\, i=1,2$. Putting $F_1(t)=u$, we get $\overline{F}_1(t)=1-u$ and $\overline{F}_2\left(Q_1(u)\right)=1-Q_3(u)$. Now quantile form of \eqref{eq3.1} is defined by
\begin{equation}
\begin{aligned}
R_\alpha^*(u) & =\int_u^1\left(\frac{1-Q_3(u)}{1-u}\right)^\alpha\left(\frac{q_1(p)}{q_1(p) q_3(p)}\right)^\alpha \frac{f_2(Q_1(p)) q_1(p)}{\left(1-Q_3(u)\right)} d p \\
& =\frac{\left(1-Q_3(u)\right)^{\alpha-1}}{(1-u)^\alpha} \int_u^1\left(q_3(p)\right)^{1-\alpha} d p.
\end{aligned}
\end{equation}
When $u \rightarrow 0, Q_3(u)=0$ so that $R_\alpha^*(u)=\int_0^1\left(q_3(p)\right)^{1-\alpha}dp=I_\alpha^*$.

\begin{rem}
    Differentiating $R_\alpha^*(u)$ with respect to $\alpha$ and substituting $\alpha=1$, we get
    \begin{equation*}
\frac{\partial R_\alpha^*(u)}{\partial \alpha}\biggr|_{\alpha=1}=\log\left(\frac{1-Q_3(u)}{1-u}\right)-\frac{1}{1-u}\int_u^1 \log q_3(p) dp,    
\end{equation*}
which is the quantile form of residual K-L divergence introduced by \cite{sankaran2016kullback}.
\end{rem}
\begin{ex}
 Suppose that $X_i$ has exponential distribution with parameter (mean) $\lambda_i,=1,2 $. Then $Q_i(u)=-\lambda_i \log (1-u), i=1,2$. Then $Q_3(u)=1-(1-u)^{\lambda_1 / \lambda_2}$ and $q_3(u)=\frac{\lambda_1}{\lambda_2}(1-u)^{\frac{\lambda_1}{\lambda_2}-1}$.
\begin{equation*}
\begin{aligned}
R_\alpha^*(u)&=\frac{(1-u)^{\frac{\lambda_1}{\lambda_2}(\alpha-1)}}{(1-u)^\alpha} \int_u^1\left(\frac{\lambda_1}{\lambda_2}\right)^\alpha(1-p)^{\left(\frac{\lambda_1}{\lambda_2}-1\right)(1-\alpha)} dp\\ 
&=(1-u)^{\frac{\lambda_1}{\lambda_2}\alpha-\frac{\lambda_1}{\lambda_2}-\alpha}\left(\frac{\lambda_1}{\lambda_2}\right)^{\alpha}\frac{(1-p)^{(\frac{\lambda_1}{\lambda_2}-1)(1-\alpha)+1}\biggr|_{u}^{1}}{-[(\frac{\lambda_1}{\lambda_2}-1)(1-\alpha)+1]}\\
& =\frac{\lambda_2\left(\lambda_1 / \lambda_2\right)^\alpha}{\lambda_2+\alpha\left(\lambda_2-\lambda_1\right)}=I_\alpha^*.
\end{aligned}
\end{equation*}   
\end{ex}

\begin{ex}
 Suppose that $X_1$ follows linear hazard quantile function with $Q_1(p)=\frac{1}{a+b} \log \left(\frac{a+bp}{a(1+p)}\right), a>0,a+b>0$ (\cite{midhu2014class}) and $X_2$ follows exponential distribution with
$Q_2(p)=-\frac{1}{\lambda_2} \log (1-p),\lambda_2>0$. Then
$Q_3(p)=1-\left(\frac{a+b p}{a(1+p)}\right)^{\frac{-\lambda_2}{a+b}}$. We now get $q_3(p)=\frac{\lambda_2}{(a+b)}(b-a)\left[\frac{(a+b p)}{a(1+p)}\right]^{\frac{-\lambda_2-a-b}{a+b}} \frac{1}{a(1+p)^2}, \quad b>a>0.$ The relative information generating function is given by
\begin{equation}
 R_\alpha^*(u)=\frac{(a+bu)^{-\frac{\lambda_2(\alpha-1)}{a+b}}}{(a(1+u))^{-\frac{\lambda_2(\alpha-1)}{a+b}}(1-u)^\alpha}\int_{u}^{1}\left(\frac{\lambda_2(b-a)}{a(a+b)}\right)^{1-\alpha}\frac{(a+bp)^{\frac{(-\lambda_2-a-b)(1-\alpha)}{a+b}}}{(a(1+p))^{\frac{(-\lambda_2+a+b)(1-\alpha)}{a+b}}}dp.
\end{equation}
There is no closed form expression for $R_\alpha^*(u)$. The above integral can be computed numerically. Letting $a=0.1,b=0.2,$ and $\lambda_2=0.6$, $R_\alpha^*(u)$ is plotted for various values of $\alpha$ and $u$ in Figure \ref{Ralphaeg}.

\begin{figure}[h!]
\centering
\includegraphics[width=7cm,height=4.6cm]{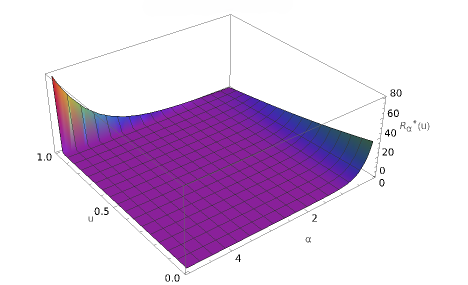}
\caption{Plot of $R_{\alpha}^*(u)$ versus $\alpha$ and $u$}
\label{Ralphaeg}
\end{figure}

\end{ex}

\begin{prop}
Let $X_1$ and $X_2$ be two non-negative random variables having $R_\alpha^*(u)$ as residual relative information generating function. A necessary and sufficient condition for $R_\alpha^*(u)=k$, a constant is that $X_1$ and $X_2$ satisfy Cox proportional hazards model property.

\begin{proof}
Assume that $R_\alpha^*(u)=k$, a constant. Then
\begin{equation}\label{3.1}
 \frac{\left(1-Q_3(u)\right)^\alpha}{(1-u)^\alpha} \frac{1}{(1-Q_3(u))} \int_u^1\left(q_3(p)\right)^{1-\alpha} dp=k.   
\end{equation}
Differentiating \eqref{3.1} with respect to $u$, we get
\begin{equation}\label{3.2}
\footnotesize
    \begin{aligned}
       \frac{\left(1-Q_3(u)\right)^{\alpha-1}}{(1-u)^\alpha}(-q_3(u))^{1-\alpha}+
       \frac{\alpha \left(1-Q_3(u)\right)^{\alpha-1}}{(1-u)^{\alpha+1}}\int_{u}^{1}q_3(p)dp-\frac{(\alpha-1)\left(1-Q_3(u)\right)^{\alpha-2}q_3(u)}{(1-u)^\alpha} \int_u^1\left(q_3(p)\right)^{1-\alpha} d p=0.
    \end{aligned}
\end{equation}
Substituting $R_\alpha^*(u)=k$ in \eqref{3.2}, we get
\begin{equation}\label{3.3}
\frac{\left(1-Q_3(u)\right)}{1-u} \alpha k-(\alpha-1) q_3(u) k=\left(\frac{1-Q_3(u)}{1-u}\right)^\alpha\left(q_3(u)\right)^{1-\alpha}    
\end{equation}
Let $\frac{q_3(u)(1-u)}{\left(1-Q_3(u)\right)}=p(u)$. Then \eqref{3.3} becomes
\begin{equation}\label{3.4}
\alpha k=k(\alpha-1) p(u)+(p(u))^{1-\alpha}.     
\end{equation}
Differentiating \eqref{3.4} with respect to $u$, we get 
\begin{equation}\label{3.6a}
    k(\alpha-1) p'(u)+(1-\alpha)(p(u))^{-\alpha}p'(u)=0,
\end{equation}
where $p'(u)$ is the derivative of $p(u)$. Now \eqref{3.6a} can be written as $(p(u))^{-\alpha}=k$ or $p(u)=k_1$. Now we can write
\begin{align}
& \frac{-d \log \left(1-Q_3(u)\right)}{d u}=-k_1 \log (1-u) \label{3.5} \\
%\Leftrightarrow \quad & 1-Q_3(u)=(1-u)^{k_1} \nonumber \\
\Leftrightarrow \quad& Q_3(u)=1-(1-u)^{k_1} \nonumber \\
\Leftrightarrow \quad & Q_2^{-1}\left(Q_1(u)\right)=1-(1-u)^{k_1} \nonumber \\
%\Leftrightarrow \quad & Q_1(u)=Q_2\left(1-(1-u)^{k_1}\right)\nonumber \\
\Leftrightarrow \quad& Q_2(u)=Q_1\left(1-(1-u)^{1/k_1}\right). \label{3.6}
\end{align}
Thus $X_1$ and $X_2$ have proportional hazards property. To prove the converse, suppose \eqref{3.6} holds.
Then, $ Q_3(u)=1-(1-u)^{k_1}$ and $ q_3(u)=k_1(1-u)^{k_1-1}$.
Thus by direct calculation, we get $R_\alpha^*(u)$ is a constant.
\end{proof}
\end{prop}

\begin{prop}
 Let $\phi_1(\cdot)$ and $\phi_2(\cdot)$ be two continuous,  non-decreasing,  and invertible functions. Then relative information generating function for $\phi_1(X_1)$ and $\phi_2(X_2)$ is obtained as 
 \begin{equation}
 \footnotesize
R_\alpha^*(\phi_1(X_1),\phi_2(X_2),u)=\frac{[1-Q_2^{-1}(\phi_2^{-1}(\phi_1(Q_1(u))))]^{\alpha-1}}{(1-u)^\alpha}\int_{u}^{1}\left[\frac{d }{dp}Q_2^{-1}(\phi_2^{-1}(\phi_1(Q_1(u))))\right]^{1-\alpha}dp
 \end{equation}
 The proof is direct.
\end{prop}

Let $G(\cdot)$ be a continuous distribution function having probability density function $g(\cdot)$ defined in the interval $[0,1]$. A general family of transformed survival models is defined by
\begin{equation}\label{gen}
\overline{F}_2(x)=G\left(\overline{F}_1(x)\right) ;\quad x>0.    
\end{equation}
Thus the baseline survival function $\overline{F}_1(\cdot)$ is transformed to the survival function $\overline{F}_2(\cdot)$ by a function $G(\cdot)$. Proportional hazards and proportional odds models are special cases of \eqref{gen}. 
The model \eqref{gen} in the quantile form is written as
$$
Q_2(p)=Q_1\left(1-G^{-1}(1-p)\right),
$$
so that
$$
Q_3(p)=Q_2^{-1}\left(Q_1(p)\right)=1-G(1-p).
$$
Thus $q_3(p)=g(1-p)$ where $g(\cdot)$ is the derivative of $G(\cdot)$ and
$$
\begin{aligned}
R_{\alpha G}^*(u)=\frac{(G(1-u))^{\alpha-1}}{(1-u)^\alpha}\int_{u}^{1}(g(1-p))^{1-\alpha} dp.
\end{aligned}
$$
In the case of proportional hazards model,  $\overline{F}_2(x)=(\overline{F}_1(x))^{\theta}$, we get  $Q_2(p)=Q_1(1-(1-p)^{1/\theta})$. Then $Q_3(p)=1-(1-p)^{\theta}$ and $q_3(p)=\theta(1-p)^{\theta-1}$. Now
\begin{equation}\label{11}
\begin{aligned}
R_{\alpha}^{*}(u)&=\frac{(1-u)^{\theta(\alpha-1)}}{(1-u)^{\alpha}}\int_{u}^{1}(\theta (1-p))^{(\theta-1)(1-\alpha)}dp\\
%& =\frac{\theta(1-u)^{(\theta-1)(1-\alpha)+1}}{(\theta-1)(1-\alpha)+1} \frac{(1-u)^{\theta(\alpha-1)}}{(1-u)^\alpha}\\
%&=\frac{\theta(1-u)^{\theta(1-\alpha)+\theta+\alpha+\theta \alpha-\theta}}{\theta(1-\alpha)+\theta+\alpha} \frac{1}{(1-u)^\alpha}\\
%&=\frac{\theta\left(1-u\right)^\theta}{\theta(1-\alpha)+\theta+\alpha}.
&=\frac{\theta^{1-\alpha}}{\theta(1-\alpha)+\alpha}.
\end{aligned}
\end{equation}
Two random variables $X_1$ and $X_2$ satisfy the proportional odds model if
\begin{equation}\label{12}
\frac{\overline{F}_2(x)}{F_2(x)} = r \frac{\overline{F}_1(x)}{F_1(x)},
\end{equation}
where $r$ is a positive constant. This is equivalent to write
\begin{equation}\label{13F}
\overline{F}_2(x) = \frac{r \overline{F}_1(x)}{1 - (1 - r) \overline{F}_1(x)}, \quad x > 0, \quad 0<r<1.
\end{equation}
The quantile version of \eqref{13F} is
\begin{equation}\label{13Q}
Q_3(p) = 1 - \frac{r (1 - p)}{1 - (1 - r)(1 - p)}, \quad 0 < p < 1.
\end{equation}
Let $G(x) = \frac{rx}{1 - (1 - \theta) x}$. Then the model reduces to the above quantile function \eqref{13Q}. We now get
\begin{equation}
    \begin{aligned}
     q_3(p) &= \frac{\left[1 - (1 - r) (1 - p)\right] r + r (1 - p) (1 - r)}{\left(1 - (1 -r) (1 - p)\right)^2}\\
&=   \frac{r}{\left(1 - (1 -r) (1 - p)\right)^2}.
    \end{aligned}
\end{equation}
Then
\[
R_{\alpha G}^*(u) = \frac{(1 - Q_3(u))^{\alpha-1}}{(1-u)^\alpha}\int_u^1 \left(\frac{r}{\left[1 - (1 - r)(1 - p)\right]^2}\right)^{1-\alpha} dp.
\]
Set $1 - (1 - r)(1 - p) = v, \quad (1 - r) dp = dv.$
Then
\begin{equation}
    \begin{aligned}
 R_{\alpha G}^{*}(u) &= \left(\frac{r(1 - u)}{1 - (1 - r)(1 - u)}\right)^{\alpha-1}\frac{r^{1-\alpha}}{(1-u)^\alpha}{\int_{1-(1 -r)(1-u)}^1 \left(\frac{1}{v^2}\right)^{1-\alpha} dv}\\
 %&=\left(\frac{r(1 - u)}{1 - (1 - r)(1 - u)}\right)^{\alpha-1}\left(\frac{r}{1-r}\right)^\alpha \frac{v^{2\alpha-1}}{2\alpha-1}\biggr|_{1-(1-r)(1-u)}^{1}\\
% &=\frac{(1-u)^{\alpha-1}}{(1-r)^{1-\alpha}}\frac{1}{(1-(1 -r)(1-u))^{\alpha-1}}\frac{1}{2\alpha-1}\left[1-[1-(1 -r)(1-u)]^{2\alpha-1}\right]
&=\frac{\left[1-[1-(1-r)(1-u)]^{2\alpha-1}\right]}{[1-(1-r)(1-u)]^{\alpha-1}(1-r)(2\alpha-1)}, \quad \alpha>1/2.
    \end{aligned}
\end{equation}

\section{Relative Information Generating Function for Past Lifetimes}\label{sec4}
 There are certain situations in lifetime analysis where past lifetimes may be taken into account rather than the residual lifetime. The past lifetime for the random variable $X_i$ is defined as $X_i^{(t)}=t-X_i \mid X_i \leq t,\quad i=1,2$. The relative generating function proposed in Section \ref{sec2} can be extended to the past lifetime random variables. The probability density function of $X_i^{(t)}$ for $i=1,2$ is given by
 \begin{equation}\label{}
 f_{X_i^{(t)}}(x)=
     \begin{cases}
     \frac{f_i(x)}{F_i(t)}& 0<x<t\\
     0 & \text{otherwise}.
     \end{cases}
 \end{equation}
Then the relative information generating function between $X_1^{(t)}$ and $X_2^{(t)}$ is defined by
\begin{equation}
 J_\alpha(x)=\int_0^t\left(\frac{f_1(x)}{F_1(t)} \biggr/\frac{f_2(x)}{F_2(t)}\right)^\alpha \frac{f_2(x)}{F_2(t)} dx.  
\end{equation}
The quantile version of $J_\alpha(x)$ is obtained as
\begin{equation}
    J_\alpha^*(u)=\frac{\left(Q_3(u)\right)^{\alpha-1}}{u^\alpha} \int_0^u\left(q_3(p)\right)^{1-\alpha} dp.
\end{equation}
When $u\rightarrow1$, $J_\alpha^*(u)\rightarrow I_\alpha^*$.

\begin{ex}
Suppose that $X_i$ follows Pareto II distribution with quantile function $Q_i(p)=(1-p)^{-\frac{1}{\beta_i}}-1, \beta_i>0, i=1,2.$.
Then $Q_3(p)=1-(1-p)^{\frac{\beta_2}{\beta_1}}$ and $q_3(p)=\frac{\beta_2}{\beta_1}(1-p)^{\frac{\beta_2}{\beta_1}-1}$.
Hence
$$
\begin{aligned}
& J_\alpha^*(u)=\frac{\left(1-(1-u)^{\frac{\beta_2}{\beta_1}}\right)^{\alpha-1}}{u^\alpha} \int_0^{u}\left(\frac{\beta_2}{\beta_1}\right)^{1-\alpha}(1-p)^{\left(\frac{\beta_2-\beta_1}{\beta_1}\right)(1-\alpha)} dp\\
& =\frac{\left(1-(1-u)^{\frac{\beta_2}{\beta_1}}\right)^{\alpha-1}}{u^\alpha}\left(\frac{\beta_2}{\beta_1}\right)^{1-\alpha} \frac{\beta_1}{\beta_2+\alpha(\beta_1-\beta_2)}\left[(1-u)^{\frac{\beta_2+\alpha(\beta_1-\beta_2)}{\beta_1}}-1\right].
\end{aligned}
$$
\end{ex}

\begin{ex}
Consider the quantile functions of Govindarajulu and reciprocal exponential distribution with quantile functions $Q_1(p)=2p-p^2$ and $Q_2(p)=-\frac{\lambda}{\log p}, \lambda>0$. Then $q_3(p)$ is obtained as
\begin{equation*}
      Q_3(p)=e^{-\frac{\lambda}{2p-p^2}}\quad \text{and}
     \quad q_3(p)=\frac{2\lambda(1-p)}{\left(2 p-p^2\right)^2}e^{-\left(\frac{\lambda}{2u-u^2}\right)}.
\end{equation*}
Now
\begin{equation*}
J_\alpha^*(u)=\frac{e^{-\frac{\lambda(\alpha-1)}{2 u-u^2}}}{u^\alpha} \int_0^u\left(\frac{2 \lambda(1-p)}{\left(2 p-p^2\right)^2}\right)^{1-\alpha}e^{-\frac{\lambda(1-\alpha)}{2 p-p^2}} d p,    
\end{equation*}
which can be evaluated numerically.  Plot of $J_{0.5}^*(u)$ for different values $\lambda$ and $u$ is shown in Figure \ref{Jalpha4.2}.

\begin{figure}[h!]
\centering
\includegraphics[width=7cm,height=5cm]{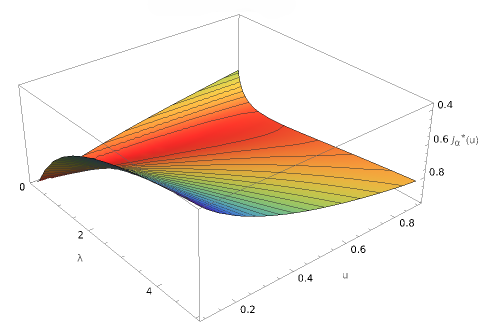}
\caption{Plot of $J_{0.5}^*(u)$ versus $\lambda$ and $u$}
\label{Jalpha4.2}
\end{figure} 
\end{ex}
   
\begin{prop}
 Let $X_1$ and $X_2$ be two non-negative random variables having $J^*_\alpha(u)$ as relative information generating function of past lifetimes. The quantity $J^*_\alpha(u) = c$, a constant, if and only if $X_1$ and $X_2$ satisfy proportional reversed hazards property.

\begin{proof}
Assume that $J^*_\alpha(u) = c$, a constant. Then we get 
\begin{equation}\label{4.4}
   (Q_3(u))^{\alpha-1} \int_0^u q_3(p)^{1-\alpha} \, dp = u^\alpha c.
\end{equation}
Differentiating \((4.4)\) with respect to $u$, we get 
\begin{equation}\label{4.5}
   \frac{\partial}{\partial u}J^*_\alpha(u)=0 \implies
(\alpha-1)(Q_3(u))^{\alpha-2}q_3(u)\int_0^u q_3(p)^{1-\alpha} \, dp + (Q_3(u))^{\alpha-1}(q_3(u))^{1-\alpha} = \alpha u^{\alpha-1} c.  
\end{equation}
Using \((4.4)\) in \((4.5)\), we obtain
\[
(\alpha - 1) u^{\alpha} c \frac{q_3(u)}{Q_3(u)} + (Q_3(u))^{\alpha-1}(q_3(u))^{1-\alpha}=\alpha u^{\alpha-1} c,
\]
which leads to
\begin{equation}\label{4.6}
    (\alpha-1)cl(u)+(l(u))^{1-\alpha}=c\alpha,
\end{equation}
 where $l(u)=\frac{uq_3(u)}{Q_3(u)}$.
 The only solution of \eqref{4.6} is $l(u)=c_1$, a constant. Thus 
 $\frac{d }{du}\log Q_3(u)=c_1\frac{d }{du}\log u$, which leads to
$ Q_3(u)=u^{c_1} \Rightarrow Q_2^{-1}(Q_1(u))=u^{c_1}$ or
$Q_1(u)=Q_2\left(u^{c_1}\right)$. 
Thus we obtain
$Q_2(u)=Q_1\left(u^{1 / c_1}\right)
$. The converse part is direct.
 \end{proof}
\end{prop}

\begin{figure}[h!]
\centering
\includegraphics[width=7cm,height=5cm]{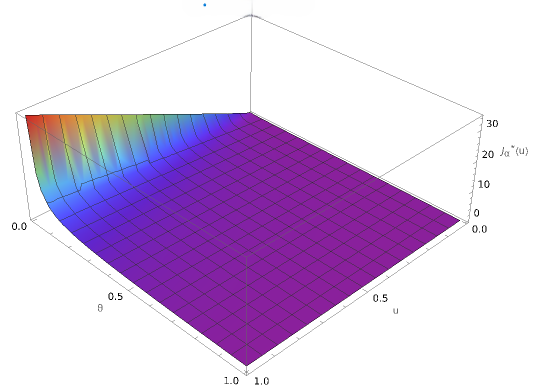}
\caption{Plot of $J_{0.7}^*(u)$ versus $\theta$ and $u$}
\label{Jalpha}
\end{figure}

Suppose that $F_2(x)=G\left(F_1(x)\right),\,\, x>0$, is a class of transformation models. Now, the quantile version of the model is
\begin{equation}\label{4.7}
Q_3(p)=Q_2^{-1}(Q_1( p))=G(p), \quad 0<p<1, 
\end{equation}
so that $q_3(p)=\frac{d G(p)}{d p}=g(p)$. Then the information generating function is
\begin{equation}
    J_{\alpha G}^*(u)=\frac{G(u)^{\alpha-1}}{u^\alpha} \int_{0}^{u}(g(p))^{1-\alpha} dp.
\end{equation}
Note that proportional reversed hazards and propotional odds models are special cases of \eqref{4.7}.
For instance, if $G(x)=\frac{x}{\theta+x(1-\theta)}$, then the model becomes proportional odds model
$
F_2(x)=\frac{F_1(x)}{\theta+F_1(x)(1-\theta)}, \,\, x>0$
which leads to
\begin{equation*}
   Q_3(u)=\frac{u}{\theta+u(1-\theta)}, \quad 0<u<1 
\end{equation*}
and 
\begin{equation*}
    q_3(u)=\frac{\theta+u(1-\theta)-u(1-\theta)}{(\theta+u(1-\theta))^2}=\frac{\theta}{(\theta+u(1-\theta))^2}.
\end{equation*}
Thus
$$
\begin{aligned}
J_\alpha^*(u)& =\frac{u ^{\alpha-1}}{(\theta+u(1-\theta))^{\alpha-1}u^\alpha} \int_0^u\left[\frac{\theta}{(\theta+p(1-\theta)}\right]^{1-\alpha} d p \\
%& \left.=\frac{\theta^{1-\alpha}}{u\left[(\theta+u(1-\theta)]^{\alpha-1}\right.} \frac{1}{1-\theta} \frac{(\theta+p(1-\theta))^\alpha}{\alpha}\right]_0^u \\
%&=\frac{\theta^{1-\alpha}}{(1-\theta) u(\theta+u(1 - \theta))^{\alpha-1}\alpha}[(\theta+u(1-\theta)-\theta] \\
%& =\frac{\theta^{1-\alpha}}{\alpha} \frac{1}{(\theta+u(1-\theta))^{\alpha-1}}.
&=\frac{1}{\alpha \theta^{\alpha}(\theta+u(1-\theta))^{\alpha-1}}\left[(\theta+u(1-\theta))^{\alpha}-\theta^\alpha\right].
\end{aligned}
$$
Plot of $J_{0.7}^*(u)$ for different values $\theta$ and $u$ is shown in Figure \ref{Jalpha}.

\section{Non parametric Estimation of Relative Information Generating Function}\label{sec5}

We now develop, non parametric estimation of the relative information generating function $R_\alpha^*(u)$ for residual lifetimes. This method can be directly applied for other two measures $I_\alpha^*$ and $J_\alpha^*(u)$. Let $Z_1, Z_2, \ldots, Z_{n}$ be independent and identically distributed samples form $Q_3(u)$. Let $Z_{(1)}<Z_{(2)}<\ldots<Z_{(n)}$ be the order statistics. Then by \cite{parzen1979nonparametric}, $Q_3(u)$ can be estimated as
\begin{equation}\label{5.1}
 \hat{Q}_3(u)=n\left(\frac{r}{n}-u\right) Z_{(r-1)}+n\left(u-\frac{r-1}{n}\right)Z_{(r)} \quad \text{for}\quad \frac{r-1}{n} \leq u \leq \frac{r}{n} \quad \text{and} \quad r=1,\ldots,n.  
\end{equation}
Then estimator of $q_3(u)$ is 
\begin{equation}\label{5.2}
 \hat{q}_3(u)=n\left(Z_{(r)}-Z_{(r-1)}\right) \quad \text{for}\quad \frac{r-1}{n} \leq u \leq \frac{r}{n} \quad; r=1,\ldots,n. 
\end{equation}
The estimator of $R_\alpha^*(u)$ is obtained as
$$
\hat{R}_\alpha^*(u)=\frac{\left(1-\hat{Q}_3(u)\right)^{\alpha-1}}{(1-u)^\alpha} \int_u^1\left(\hat{q}_3(p)\right)^{1-\alpha} dp.
$$
Thus the non parametric estimator of $R_\alpha^*(u)$ is
\begin{equation}\label{5.3}
   \hat{R}_\alpha^*(u)=\frac{\left(1-\hat{Q}_3(u)\right)^{\alpha-1}}{(1-u)^\alpha} \frac{1}{n}\left[\sum_{j:Z_{(j-1)}\leq u <Z_{(j)}}^{n}[n(Z_{(j)}-Z_{(j-1)})]^{1-\alpha}\right]. 
\end{equation}
Specifically with $u=0$, the non parametric estimator of $I_\alpha^*$ is obtained as
\begin{equation}\label{5.4}
   \hat{I}_\alpha^*=\frac{1}{n}\left[\sum_{1}^{n}[n(Z_{(j)}-Z_{(j-1)})]^{1-\alpha}\right]. 
\end{equation}
On similar lines, $J_\alpha^*(u)$ can be estimated as
\begin{equation}\label{5.5}
\begin{aligned}
    \hat{J}_\alpha^{*}(u)&=\frac{\left(\hat{Q}_3(u)\right)^{\alpha -1}}{u^ \alpha} \int_0^{u}\left(\hat{q}_3(p)\right)^{1-\alpha} dp\\
    &=\frac{\left(\hat{Q}_3(u)\right)^{\alpha -1}}{u^ \alpha} \frac{1}{n}\left[\sum_{1}^{j:Z_{(j-1)}\leq u <Z_{(j)}}[n(Z_{(j)}-Z_{(j-1)})]^{1-\alpha}\right].
\end{aligned}
\end{equation}

\section{Simulation Study}\label{sec6}
A simulation study is conducted to assess the performance of the estimator of the relative information generating function. We consider the estimators of $I_\alpha^*$ and $R_\alpha^*(u)$. The same procedure can be applied for $J_\alpha^*(u)$. Suppose that $X_1$ follows Govindarajulu distribution with quantile function $Q_1(p) = 2p-p^2$ and $X_2$ follows a reciprocal exponential distribution with $Q_2(p) = -\frac{-\lambda}{\log p}$. Then $Q_3(p) = e^{-\frac{\lambda}{2p-p^2}}$ and $q_3(p) = e^{-\frac{\lambda}{2p-p^2}}\frac{2\lambda(1-p)}{(2p-p^2)^2}$.
Then
\begin{equation*}
    I_\alpha^*=\int_{0}^{1}(2\lambda(1-p))^{1-\alpha}(2p-p^2)^{-2(1-\alpha)}e^{\frac{-\lambda(1-\alpha)}{2p-p^2}}dp
\end{equation*}
and
\begin{equation*}
    R_\alpha^*(u)=\frac{(1-e^{\frac{-\lambda}{2u-u^2}})}{(1-u)^\alpha}\int_{u}^{1}(2\lambda(1-p))^{1-\alpha}(2p-p^2)^{2(1-\alpha)}e^{\frac{-\lambda(1-\alpha)}{2p-p^2}}dp.
\end{equation*}
We take $\lambda = 0.7$ and $\alpha = 0.3$. The proposed estimators of $I_\alpha^*$ and $R_\alpha^*(u)$ are calculated using 1000 Monte Carlo simulations. 

The bias and mean squared error (MSE) of the estimators of $I_\alpha^*$ for $n = 50, 100, 250$, and $500$ are calculated and presented in Table \ref{I}. The same for the estimator of $R_\alpha^*(u)$ for various values of $u = 0.25, 0.5$, and $0.75$, with sample size $n = 50,100,250,$ and $500$ are presented in Table \ref{Ralphasim}. From the tables, we observe that the MSE decreases with increasing sample sizes and the estimators are asymptotically unbiased. As the performance of the estimators of $I_\alpha^*$ and $R_\alpha^*(u)$ for other values of $\alpha$ is similar, corresponding results are not included here. Similar results are obtained for $J_\alpha^*(u)$ as well.

% Please add the following required packages to your document preamble:
% \usepackage{multirow}
% Please add the following required packages to your document preamble:
% \usepackage{multirow}

\begin{table}[h!]
\centering
\caption{Simulation results for Information generating function $I_\alpha^*$}
%\label{frequent_fixed_scenario2}
\label{I}
\begin{tabular}{|c|c|c|}
\hline
n   & Bias   & MSE    \\ \hline
50  & 0.3436 & 0.1184 \\ 
100 & 0.3424 & 0.1174 \\ 
250 & 0.3423 & 0.1172 \\ 
500 & 0.3418 & 0.1168 \\ \hline
\end{tabular}
\end{table}

\begin{table}[h!]
\centering
\caption{Simulation results for $R_\alpha^*(u)$}
\label{Ralphasim}
\begin{tabular}{|c|cc|cc|cc|}
\hline
\multirow{2}{*}{n} & \multicolumn{2}{c|}{u=0.25}          & \multicolumn{2}{c|}{u=0.5}           & \multicolumn{2}{c|}{u=0.75}          \\ \cline{2-7} 
                   & \multicolumn{1}{c|}{Bias}   & MSE    & \multicolumn{1}{c|}{Bias}   & MSE    & \multicolumn{1}{c|}{Bias}   & MSE    \\ \hline
50                 & \multicolumn{1}{c|}{0.2586} & 0.0687 & \multicolumn{1}{c|}{0.0847} & 0.0098 & \multicolumn{1}{c|}{0.0318} & 0.0024 \\ 
100                & \multicolumn{1}{c|}{0.2303} & 0.0540 & \multicolumn{1}{c|}{0.0687} & 0.0060 & \multicolumn{1}{c|}{0.0109} & 0.0007 \\ 
250                & \multicolumn{1}{c|}{0.2238} & 0.0505 & \multicolumn{1}{c|}{0.0598} & 0.0041 & \multicolumn{1}{c|}{0.0051} & 0.0003 \\ 
500                & \multicolumn{1}{c|}{0.2179} & 0.0477 & \multicolumn{1}{c|}{0.0562} & 0.0034 & \multicolumn{1}{c|}{0.0007} & 0.0001 \\ \hline
\end{tabular}
\end{table}

\section{Application to Prostate Cancer Data}\label{sec7}
This study is based on data from a randomized clinical trial assessing treatments for stage 3 and stage 4 prostate cancer. Patients received placebo pills or diethylstilbestrol (DES) at doses of 0.2 mg, one mg, or five mg. We focus on the survival times of individuals treated with placebo ($X_1$) and those treated with five mg of DES ($X_2$), as per \cite{andrews2012data}.

\begin{figure}[h!]
\centering
\includegraphics[width=7cm,height=5cm]{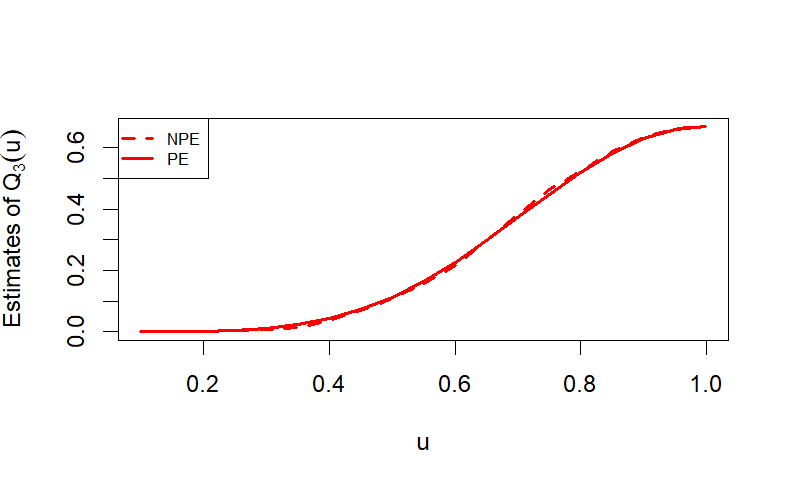}
\caption{The estimates of $Q_{3}(u)$ versus $u$}
\label{Q3}
\end{figure}

\begin{figure}[h!]
\centering
\includegraphics[width=7cm,height=5cm]{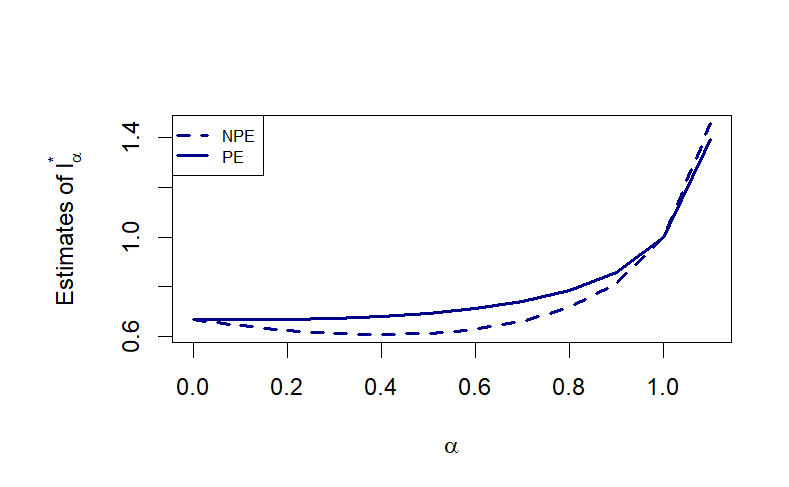}
\caption{The estimates of $I_{\alpha}^*$ versus $\alpha$}
\label{Ialphadata}
\end{figure}

\cite{un2022some} have analysed the dataset and fitted parametric models to the data using the method of $L$-moments. Govindarajulu distribution (\cite{nair2012govindarajulu} and power distributions were identified as appropriate fits to variables $X_1$ and $X_2$ respectively, having the estimated quantile functions 
\begin{equation*}
     Q_1(u)=69.26 (2.04\,u^{1.04}-1.04\,u^{2.04}), \quad 0\leq u\leq 1,
\end{equation*}
and 
\begin{equation*}
     Q_2(u)=74.13 \, u^{0.17}, \quad 0\leq u\leq 1.
\end{equation*}
Therefore, 
\begin{equation}\label{dataQ3}
    Q_3(u)=Q_2^{-1}(Q_1(u))=\left[ \frac{69.26 \left\{ \left( 2.04\right)u^1.04-1.04 u^{2.04} \right\}}{74.13} \right]^{\frac{1}{0.17}}.
\end{equation} 
This parametric model for $Q_3(u)$ is used to generate the random sample $Z_1, \dots, Z_n$ for subsequent analysis. In our first study, the parametric estimate (PE) of  $Q_3(u)$ in \eqref{dataQ3} is compared with the corresponding non parametric estimate (NPE) computed by \eqref{5.1}. The comparison is shown in Figure \ref{Q3}. Furthermore, on differentiating \eqref{dataQ3} with respect to $u$, $q_3(u)$ is derived and substituting in \eqref{2.2}, PE of $I_{\alpha}^*$ is obtained. Its NPE given by \eqref{5.4} is also obtained and shown in Figure \ref{Ialphadata}. As depicted in Figures \ref{Q3} and \ref{Ialphadata}, the plots suggest strong alignment between NPE and PE, demonstrating the effectiveness of the proposed non parametric methods.

% \begin{figure}[h!]
% \centering
% \includegraphics[width=7cm,height=4cm]{qd3.png}
% \caption{The estimates of $q_{3}(u)$ versus $u$}
% \label{q3}
% \end{figure}

\begin{figure}[h!]
\centering
\includegraphics[width=7cm,height=4cm]{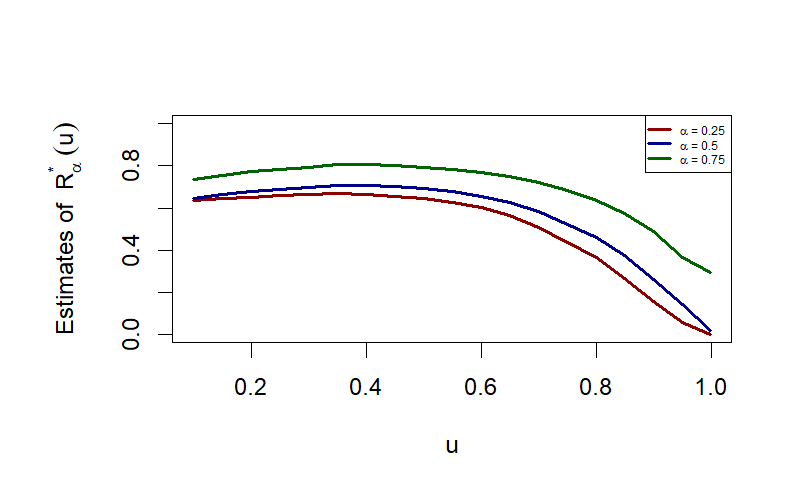}
\caption{The estimates of $R_{\alpha}^{*}(u)$ versus $u$}
\label{Ralphadata}
\end{figure}

\begin{figure}[h!]
\centering
\includegraphics[width=7cm,height=4cm]{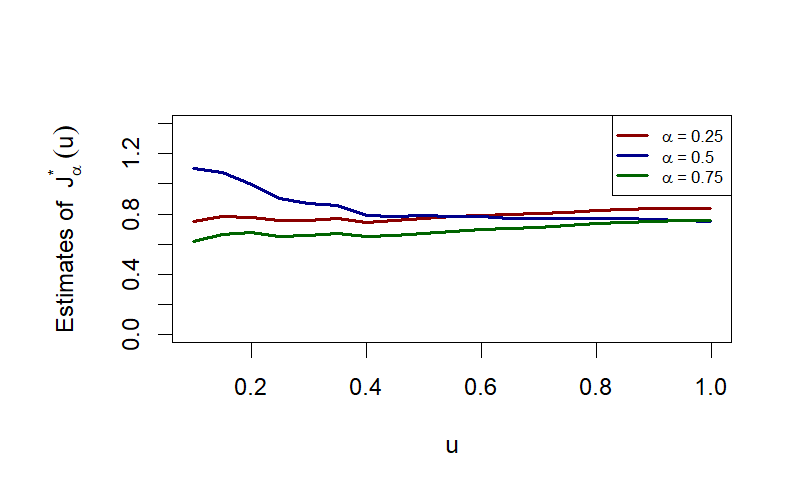}
\caption{The estimates of $J_{\alpha}^{*}(u)$ versus $u$}
\label{Jalphadata}
\end{figure}

Additionally, the NPE of $R_{\alpha}^*(u)$ given by  \eqref{5.3} and $J_{\alpha}^*(u)$ given by \eqref{5.5} are plotted against $u$ for $\alpha=0.25,0.5,$ and $0.75$ in Figures 
\ref{Ralphadata} and \ref{Jalphadata} respectively.


The relationships between $I_{\alpha}^*$ and various divergence measures, established in Section \ref{sec2}, are used to non parametrically estimate some of these measures. The results are presented in Table \ref{data_DM}, which also includes divergence estimates for survival times of individuals receiving placebo and one mg dose of DES. These estimates indicate that higher DES dose increases the disparity between placebo and DES, providing strong evidence to distinguish the placebo from the five mg dose. Furthermore, the findings highlight that dosage variations significantly impact survival times.

\begin{table}[]
\centering
\caption{Estimated divergence measures}
\label{data_DM}
\begin{tabular}{|c|c|c|}
\hline
                       & Placebo v/s 5 mg      & Placebo v/s 1 mg      \\ \hline
K-L  divergence                  & 2.6851                & 1.2256                \\ \hline
Hellinger distance                  & 0.3866                & 0.2574               \\ \hline
Bhattacharyya distance                   & 0.4888               &     0.2976            \\ \hline
R$\acute{e}$nyi measure of order 0.25                & 0.5329                 &   0.2860           \\ \hline
R$\acute{e}$nyi measure of order 0.5                  & 0.7994                &  0.4290              \\ \hline
R$\acute{e}$nyi measure of order 0.75                 & 1.5989               & 0.8581                \\ 
 \hline
\end{tabular}
\end{table}

\section{Conclusion}
The paper proposed a quantile based relative information generating function and its dynamic versions in residual and past lifetimes. The proposed function provides a method to derive divergence measures between two probability distributions, when there is tractable form for either of distribution functions. Several properties of the function were derived. Various examples were presented to illustrate the importance of the measure. Non parametric estimator of the relative information generating function was derived using the quantile density estimator of \cite{parzen1979nonparametric}. The efficiency of the estimator was studied through a simulation work. The proposed methodology was then applied to a real data on prostate cancer.

\section*{Acknowledgement}

The third author wishes to acknowledge the financial support of the Council of Scientific \& Industrial Research, Government of India,  via the Junior Research Fellowship scheme under reference No. 09/0239(13499)/2022-EMR-I.

\section*{Disclosure statement}
On behalf of all authors, the corresponding author states that there are no conflicts of interest between the authors.

\bibliographystyle{apalike}
\bibliography{ref}
\end{document}